# The frequency problem of the three gap theorem

HUIXING ZHANG[1]

**ABSTRACT:** The three gap theorem was originally a conjecture by Steinhaus, who asserted that there are at most three distinct gap lengths in the fractional parts of the sequence $\alpha, 2\alpha, \cdots N\alpha$, for any integer $N$ and real number $\alpha$. This conjecture has been proved by many different methods, and extended to higher dimensional cases. For the three gap theorem, what is the frequence of two gaps or three gaps? On the basis of Ravenstein's[1] work, this paper studies the frequency of the three gap theorem through cake model. We show that for almost all $\alpha \in (0,1) \setminus \mathbb{Q}$ in the sence of the Lebesgue measure, the frequency of the two gaps is 0, in other words, the frequency of the three gaps is 1. The proof makes full use of the continued fraction expansion of $\alpha$, as well as partial results in the Diophantine approximation. Finally, we illustrate that two gaps occur with a frequency of 0.

## 0  Introduction

The three-gap theorem was originally a conjecture made by Steinhaus, who asserted that for any integer $N$ and real number $\alpha$, the fractional part of the sequence $\alpha$, $2\alpha$, $\cdots$, $N\alpha$, that is,

$$\langle\alpha\rangle, \langle 2\alpha\rangle, \cdots, \langle N\alpha\rangle$$

there are at most three different gap lengths. A more refined result is that if $\alpha$ is irrational, the above sequence has either two different gap lengths or three different gap lengths. Throughout the paper, the symbol $\langle x \rangle$ will be fixed to denote the fractional part of $x$, that is, $\langle x \rangle = x \bmod 1$.

Since Steinhaus formulated the conjecture, many scholars have proved the three gap theorem using different methods [2][3]. In 1958, Świerckowski[4] proved a part of Steinhaus's conjecture by means of a recurrence relation. In 1967, Slater[5] formulated the step problem and the gap problem, and the corresponding proof is given by using continued fractions. In 1988, Ravenstein[1] proved the three gap theorem using the recurrence relation given by Świerckowski, and established a connection between the results of his proof and the continued fraction, which was used to carve out the two dividing points closest to the origin.

In 1988, Ravenstein[1] published a paper entitled 《The three gap theorem (Steinhaus conjecture)》, in which the author equated the problem of partitioning the partition of the partition point $\{\langle\alpha\rangle,\langle 2\alpha\rangle,\cdots,\langle n\alpha\rangle\}$ with respect to the interval $[0,1]$ to the problem of distributing the distribution of $n$ consecutively placed points on the unit circumference with $\alpha$ as the angle of rotation, and thus inscribed the cake model of the Three Gap Theorem, i.e., giving another geometric interpretation of the Three Gap Theorem. Using the cake model, we will investigate the frequency problem of the three-gap theorem.

It may be assumed that n=13 in the partition $\{\langle\alpha\rangle,\langle 2\alpha\rangle,\cdots,\langle n\alpha\rangle\}$ and that there exists $\alpha \in (0,1)\setminus\mathbb{Q}$, such that the first circle is cut with at most five cuts, then the corresponding cake model diagram is shown below:

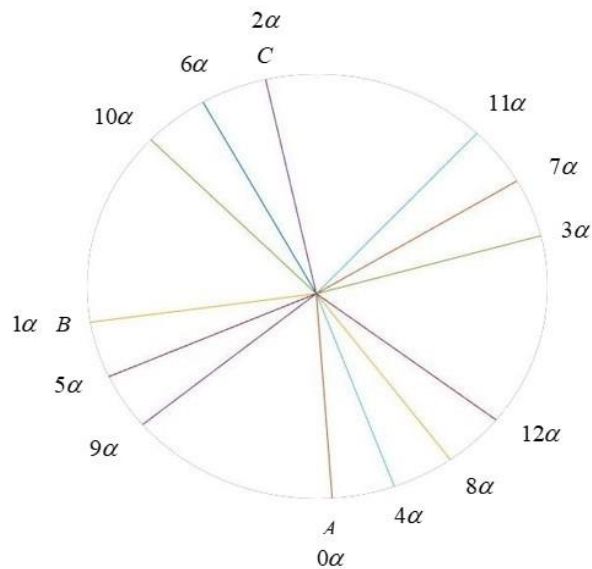

Fig 1: Cake model

Think of the unit circle as a complete cake. Cut the first knife at point A of the unit circumference, the cutting length is the radius of the unit circumference, point A is $0\alpha$; clockwise rotation angle $\alpha$ cut the second knife, and the unit circumference intersected at point B, point B is $1\alpha$; and then clockwise rotation angle $\alpha$ cut the third knife, point C is $2\alpha$; and then in turn, clockwise rotation angle $\alpha$, to get the unit circumference of the point $3\alpha$, $4\alpha$; This is the end of the cutting of the first circle. Next start the second lap of the cut with $4\alpha$. At point $4\alpha$ clockwise rotation angle $\alpha$ after cutting a knife, get the unit circle on the point $5\alpha$; and then clockwise rotation angle $\alpha$, get the unit circle on the point $6\alpha$, $7\alpha$, $\cdots$, $12\alpha$. So when $n=13$, can cut the unit circumference of 13 cuts, the last cut point for $12\alpha$.

Consider $\langle k\alpha\rangle (0 \le k \le 12)$ to be the arc length distance between the point $k\alpha$ and the starting point $0\alpha$ (the direction of the arc length distance is counterclockwise).

Therefore, for the cake model, the gap length in the three gap theorem refers to the lengths of those arc segments of the unit circumference that are divided by the sequence $\{\langle k\alpha \rangle\}_{k=0}^{n-1}$. Therefore, from the conclusion of the three gap theorem, there are only two or three cases of the size of these arc segment lengths.

# 1 Description of the problem

In order to be able to have a deep understanding of the research in this chapter, a clear description of the problem under study is given here first.

Let $\alpha \in (0,1) \backslash \mathbb{Q}$, define the sequence $P_N(\alpha)$ as

$$P_N(\alpha) = \{\langle 0\alpha \rangle, \langle 1\alpha \rangle, \langle 2\alpha \rangle, \cdots, \langle (N-1)\alpha \rangle, 1\}, \text{ where } N \in \mathbb{N}.$$

after arranging the points in the sequence $P_N(\alpha)$ in ascending order, denote the sequence of sets with the length of the gap between two points as $\Delta_N(\alpha)$, Let

$$\Delta_N(\alpha) = \{m_1, m_2, \cdots\cdots, m_N\}.$$

**Note 1.** Set $\{2,2\}$ and set $\{2\}$ are the same set.

Noting that the number of distinct elements in $\Delta_N(\alpha)$ is $\#\Delta_N(\alpha)$, it follows from the three gap theorem that

$$\#\Delta_N(\alpha) \in \{2,3\}.$$

That is, there are either two or three gap lengths formed after the sequence $P_N(\alpha)$ is partitioned into the interval $[0,1]$. Then a natural question will arise: how often are there two gap lengths and three gap lengths in the process of dividing the interval $[0,1]$ by the sequence $P_N(\alpha)$?

It can be known that when $N \to \infty$, the frequency of occurrence of both gaps is

$$\lim_{N \to \infty} \frac{\#\{1 \leq n \leq N \mid \#\Delta_n(\alpha) = 2\}}{N}.$$

in fact, it appears that once the frequency results for the two gaps are calculated, the frequency results for the three gaps will naturally be derived as well.

In response to the above question, a result is derived for the frequency of occurrence of the two gaps almost everywhere in the sense of the Lebesgue measure, i.e., the following theorem.

**Theorem 1.** For almost all $\alpha \in (0,1) \backslash \mathbb{Q}$ in the sense of the Lebesgue measure, there are

$$\lim_{N\to\infty}\frac{\#\{1\leq n\leq N\mid \#\Delta_n(\alpha)=2\}}{N}=0.$$

## 2  Background knowledge

A proof of the three gap theorem was given in 1988 in a paper by Ravenstein[1] entitled 《The three gap theorem (Steinhaus conjecture)》. The authors equate the problem of partitioning a sequence $P_N(\alpha)$ into the interval $[0,1]$ to the problem of distributing $N$ points placed consecutively on the unit circumference with $\alpha$ as the angle of rotation, thus giving the cake model and an interesting geometric interpretation of the simple continuous fractional spread of $\alpha$. The authors use a labeled sequence to further give a proof of the three gap theorem by means of a recursive relation[4], and give a relation between two of the division points and a simple connected fractional expansions of $\alpha$.

Define a symbol $\{u_j(N)\}_{j=1}^{N}$ to be a family of increasing sequences satisfying

$$\langle u_j\alpha\rangle < \langle u_{j+1}\alpha\rangle, \quad \text{where} \quad j=1,2,\cdots,N.$$

generally, $u_j(N)$ is abbreviated as $u_j$. That is, after arranging the points in the sequence $P_N(\alpha)$ in ascending order as

$$\begin{aligned}P_N(\alpha)&=\{\langle 0\alpha\rangle,\langle 1\alpha\rangle,\langle 2\alpha\rangle,\cdots,\langle(N-1)\alpha\rangle,1\}\\&=\{\langle u_1\alpha\rangle,\langle u_2\alpha\rangle,\langle u_3\alpha\rangle,\cdots,\langle u_N\alpha\rangle,1\}\end{aligned}.$$

Corresponding to the cake model Figure 1, there is

$$u_1=0,u_2=9,u_3=5,u_4=1,\cdots,u_{12}=8,u_N=u_{13}=4.$$

When cutting with $A$ (i.e., the point $0\alpha$), the point $u_2\alpha$ and the point $u_N\alpha$ are both points close to $A$. In the work of Ravenstein[1], it is the relation between $u_2$、 $u_N$ and the simple continued fractional expansions of $\alpha$ that is given.

### 2.1  Simple continued fractional expansions of $\alpha$

In this section, some knowledge about simple continued fractional expansions of $\alpha$ will be presented[6].

Call a function

$$a_0 + \cfrac{1}{a_1 + \cfrac{1}{a_2 + \cfrac{1}{a_3 + \cdots + \cfrac{1}{a_n}}}}$$

of $n+1$ variables

$$a_0, a_1, a_2, \cdots\cdots, a_n$$

a finite continued ftaction. And when $n \to \infty$, the above function is said to be an infinitely connected fraction. Continued fractions are important in many branches of mathematics, for finite continued fractions, especially in the theory of approximation of the real numbers by rational numbers; the main significance for infinite continued fractions lies in their application to the representation of irrational numbers.

The formula for continued fractions is cumbersome and inconvenient, and we usually use

$$[a_0; a_1, a_2, a_3, \cdots, a_n]$$

to remember finite connected fractions. Call $a_1, a_2, a_3, \cdots\cdots, a_n$ the partial quotient of a continued fraction, or simply the quotient. It is common to use

$$[a_0; a_1, a_2, a_3, \cdots]$$

to memorize the infinite series of scores. Call

$$a_n' = [a_n, a_{n+1}, \cdots]$$

the nth perfect quotient of a continued fraction $\alpha = [a_0; a_1, a_2, a_3, \cdots]$. Thus it can be written as

$$\alpha = [a_0; a_1, a_2, a_3, \cdots, a_{n-1}, a_n'].$$

A cyclic continued fraction is an infinite continued fraction in which

$$a_l = a_{l+k}$$

for some fixed positive number $k$ and for all $l \geq L$.

The set consisting of the partial quotients

$$a_L, a_{L+1}, \cdots\cdots, a_{L+k-1}$$

is called a cycle, and the cycle continued fraction can be written as

$$[a_0; a_1, \cdots\cdots, a_{L-1}, \dot{a}_L, a_{L+1}, \cdots\cdots, \dot{a}_{L+k-1}].$$

If

$$a_0 \in \mathbb{Z}; \quad a_1, a_2, \cdots\cdots, a_n \in \mathbb{Z}_+,$$

then a continued fraction is said to be a simple continued fraction.

In fact, every irrational number can be expressed in the only way possible as an infinitely simple continued fraction; and the value of an infinitely simple continued fraction must be an irrational number.

Thus, a simple continued fractional expansion of $\alpha \in (0,1) \backslash \mathbb{Q}$ can be expressed as

$$\alpha = a_0 + \cfrac{1}{a_1 + \cfrac{1}{a_2 + \cfrac{1}{a_3 + \cdots}}} = [a_0; a_1, a_2, a_3, \cdots]$$

where

$$a_0 = 0; \quad a_1, a_2, a_3, \cdots \in \mathbb{Z}_+.$$

These partial quotients $a_i$ satisfy the following recurrence relation:

$$\begin{cases} p_n = a_n p_{n-1} + p_{n-2} \\ q_n = a_n q_{n-1} + q_{n-2} \end{cases},$$

where

$$\begin{cases} p_{-2} = 0 \\ q_{-2} = 1 \end{cases}, \begin{cases} p_{-1} = 1 \\ q_{-1} = 0 \end{cases}, \begin{cases} p_0 = a_0 \\ q_0 = 1 \end{cases}, \begin{cases} p_1 = a_0 a_1 + 1 \\ q_1 = a_1 \end{cases}.$$

Call

$$\frac{p_n}{q_n} = [a_0; a_1, a_2, a_3, \cdots, a_n]$$

is the nth asymptotic fraction of $\alpha$. Moreover, $p_n$ and $q_n$ can be represented as matrices:

$$\begin{pmatrix} p_{n+1} & p_n \\ q_{n+1} & q_n \end{pmatrix}$$
$$= \begin{pmatrix} p_n & p_{n-1} \\ q_n & q_{n-1} \end{pmatrix} \begin{pmatrix} a_{n+1} & 1 \\ 1 & 0 \end{pmatrix}$$
$$= \begin{pmatrix} a_0 & 1 \\ 1 & 0 \end{pmatrix} \begin{pmatrix} a_1 & 1 \\ 1 & 0 \end{pmatrix} \begin{pmatrix} a_2 & 1 \\ 1 & 0 \end{pmatrix} \cdots \begin{pmatrix} a_n & 1 \\ 1 & 0 \end{pmatrix} \begin{pmatrix} a_{n+1} & 1 \\ 1 & 0 \end{pmatrix}.$$

In fact, the asymptotic fraction $\dfrac{p_n}{q_n}$ of $\alpha$ converges exactly to $\alpha$, where

$$\frac{p_n}{q_n} = \frac{p_{n,a_n}}{q_{n,a_n}} = \frac{a_n p_{n-1} + p_{n-2}}{a_n q_{n-1} + q_{n-2}}.$$

We say that $\dfrac{p_{n,i}}{q_{n,i}}$ locally converges to $\alpha$, if

$$\frac{p_{n,i}}{q_{n,i}} = \frac{ip_{n-1} + p_{n-2}}{iq_{n-1} + q_{n-2}} = [a_0; a_1, a_2, \cdots\cdots, a_{n-1}, i], \quad i = 1, 2, \cdots\cdots, a_n;$$

where

$$p_{-2} = q_{-1} = 0, \quad q_{-2} = p_{-1} = 1.$$

### 2.2 Two important conclusions

In order to solve the frequency problem of the three gap theorem, we need to first introduce two conclusions from the work of Ravenstein[1]. The first conclusion is the second part of the literature [1], i.e. Lemma 2, which gives a recursive relation between the division points by means of the ordered sequence $\{u_j\alpha\}_{j=1}^{N}$.

**Lemma 2.** If $N = u_2 + u_N$, then $u_j = ((j-1)u_2) \bmod N$, where $j = 1, 2, \cdots, N$.

Lemma 2 implies that when a $N = u_2 + u_N$, there is $\#\Delta_N(\alpha) = 2$. In other words, at this point the length of the gap formed by the division of the unit circle by the sequence $P_N(\alpha)$ will have two different sizes.

The second conclusion is the third part of the paper, i.e. the following Theorem 2, which gives the values of $u_2$ and $u_N$ from a simple continued fraction spread of $\alpha$.

**Theorem 3.**

(i) If $q_{n,i-1} < N \leq q_{n,i}$, where $2 \leq i \leq a_n$ ($n \geq 2$), then

$$u_2 = \begin{cases} q_{n-1}, \text{若n为奇数} \\ q_{n,i-1}, \text{若n为偶数} \end{cases}$$

$$u_N = \begin{cases} q_{n,i-1}, \text{若n为奇数} \\ q_{n-1}, \text{若n为偶数} \end{cases}.$$

(ii) If $q_{n-1} < N \leq q_{n,1}$ ($n \geq 2$), then

$$u_2 = \begin{cases} q_{n-1}, \text{若n为奇数} \\ q_{n-2}, \text{若n为偶数} \end{cases}$$

$$u_N = \begin{cases} q_{n-2}, \text{若n为奇数} \\ q_{n-1}, \text{若n为偶数} \end{cases}$$

(iii) If $N \leq q_1$, then $u_j = j - 1$, $j = 1, 2, \cdots\cdots, N$.

With respect to Theorem 3, it can be understood using the following line diagram:

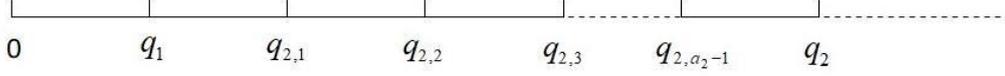

Fig 2: The values of $u_2$ and $u_N$

When $n=1$: if $N \leq q_1$, then $u_2 = 1$, $u_N = N-1$.

When $n=2$: if $q_1 < N \leq q_{2,1}$, then $u_2 = q_0$, $u_N = q_1$; if $q_{2,1} < N \leq q_{2,2}$, then $u_2 = q_{2,1}$, $u_N = q_1$; if $q_{2,2} < N \leq q_{2,3}$, then $u_2 = q_{2,2}$, $u_N = q_1$; $\cdots$ ; if $q_{2,a_2-1} < N \leq q_2$, then $u_2 = q_{2,a_2-1}$, $u_N = q_1$.

When $n=3$ 时: if $q_2 < N \leq q_{3,1}$, then $u_2 = q_2$, $u_N = q_1$; if $q_{3,1} < N \leq q_{3,2}$, then $u_2 = q_2$, $u_N = q_{3,1}$; $\cdots$

……

In other words, when $\alpha$ is given, $a_0, a_1, a_2, \cdots, a_n$ will be given, and then $q_1, q_{2,1}, q_{2,2}, q_{2,3}, \cdots$ will be given. By Theorem 3, if $N$ is in the interval where $u_2$ and $u_N$ are given, then $u_2 + u_N$ is also given. So it is only necessary to find out which $N$ in the interval where $N$ is located will make the conclusion that $N = u_2 + u_N$ hold, and then the size of the gap length at the Nth cut at this point in time will be in two cases.

### 2.3 Proof of Theorem 1

In this section, we give a detailed proof of Theorem 1. See the literature[7] for the following important results.

**Lemma 4.** For almost all real numbers $x = [a_1, a_2, \cdots, a_n] \in (0,1)$, there are

$$\lim_{n \to \infty} \frac{\log q_n(x)}{n} = \frac{\pi^2}{12 \log 2}.$$

Therefore, if $x$ satisfies

$$\lim_{n \to \infty} \frac{\log q_n(x)}{n} = \frac{\pi^2}{12 \log 2} =: t.$$

then

$$\lim_{n\to\infty}\frac{\log q_n(x)}{nt}=1.$$

thus

$$\log q_n(x)\sim nt\,(n\to\infty),\quad i.e.\ q_n(x)\sim 10^{nt}\,(n\to\infty).$$

Here is the Borel-Bernstein theorem, see literature [8]

**Lemma 5.**[8] （Borel-Bernstein）Let $\Phi$ be an arbitrary positive function defined on the set of natural numbers $\mathbb{N}$, define

$$E(\Phi)=\{x\in[0,1):a_n(x)\geq\Phi(n)\ i.o.\ n\},$$

then

$$Leb(E(\Phi))=\begin{cases}0,\text{如果}\sum_{n=1}^{\infty}\dfrac{1}{\Phi(n)}\text{收敛}\\[2mm]1,\ \text{如果}\sum_{n=1}^{\infty}\dfrac{1}{\Phi(n)}\text{不收敛}\end{cases}.$$

**Lemma 6.** For almost all $\alpha\in(0,1)\setminus\mathbb{Q}$ in the sense of the Lebesgue measure, there are

$$\lim_{n\to\infty}\frac{a_1+a_2+\cdots+a_n}{q_{n-1}}=0.$$

where $q_{n-1}(\alpha)$ is abbreviated to $q_{n-1}$.

**Proof.** By lemma 4,

$$\lim_{n\to\infty}\frac{a_1+a_2+\cdots+a_n}{q_{n-1}}=\lim_{n\to\infty}\frac{a_1+a_2+\cdots+a_n}{10^{(n-1)t}}.$$

By lemma 5, let

$$(E(\Phi))^c=\{\alpha\in(0,1)\setminus\mathbb{Q}:\exists N\in\mathbb{N},\text{当}n\geq N\text{时},\ \text{有}a_n(\alpha)<\Phi(n)\}.$$

then

$$Leb\big((E(\Phi))^c\big)=\begin{cases}1,\text{如果}\sum_{n=1}^{\infty}\dfrac{1}{\Phi(n)}\text{收敛}\\[2mm]0,\ \text{如果}\sum_{n=1}^{\infty}\dfrac{1}{\Phi(n)}\text{不收敛}\end{cases}.$$

Since $\Phi$ is an arbitrary positive function, it is useful to take $\Phi(n)=n^2$, then

$$\sum_{n=1}^{\infty}\frac{1}{\Phi(n)}=\sum_{n=1}^{\infty}\frac{1}{n^2}$$

is convergent. Thus
$$Leb\left(\left(E(\Phi)\right)^c\right) = 1.$$

Therefore
$$\lim_{n \to \infty} \frac{a_1 + a_2 + \cdots + a_n}{q_{n-1}}$$
$$= \lim_{n \to \infty} \frac{a_1 + a_2 + \cdots + a_{N-1} + a_N + a_{N+1} + \cdots + a_n}{q_{n-1}}.$$

Let
$$a_1 + a_2 + \cdots + a_{N-1} = v,$$

then
$$\lim_{n \to \infty} \frac{a_1 + a_2 + \cdots + a_n}{q_{n-1}} = \lim_{n \to \infty} \frac{v + a_N + a_{N+1} + \cdots + a_n}{10^{(n-1)t}}$$
$$< \lim_{n \to \infty} \frac{v + n^2 \cdot (n - N + 1)}{10^{(n-1)t}}$$
$$= 0.$$

So we prove Lemma 6. □

Next, we will prove Theorem 1

**Proof.** By theorem 3:

Scenario 1: $q_{n,i-1} < N \leq q_{n,i}$, where $2 \leq i \leq a_n$ and $n \geq 2$。

Whether n is odd or even, there is
$$u_2 + u_N = q_{n-1} + q_{n,i-1} = q_{n-1} + (i-1)q_{n-1} + q_{n-2}$$
$$= iq_{n-1} + q_{n-2} = q_{n,i}.$$

So, given $\alpha$ and $n$, it follows from $2 \leq i \leq a_n$ that
$$\left\{(q_{n,i-1}, q_{n,i}]\right\}_{i=2}^{a_n} = \left\{(q_{n,1}, q_{n,2}], (q_{n,2}, q_{n,3}], \cdots, (q_{n,a_n-1}, q_{n,a_n}]\right\}.$$

Therefore, for any interval $R \in \left\{(q_{n,i-1}, q_{n,i}]\right\}_{i=2}^{a_n}$, there exists a unique $N \in R$ such that $N = u_2 + u_N$. In other words, in each small interval in the set $\left\{(q_{n,i-1}, q_{n,i}]\right\}_{i=2}^{a_n}$, there is and only one $N$, such that the cut to the Nth cut is of two gaps.

So, given $\alpha$ and $n$, when $q_{n,i-1} < N \leq q_{n,i}$ ($2 \leq i \leq a_n$), there are at most $(a_n - 1)$ times when two gaps occur.

Scenario 2: $q_{n-1} < N \leq q_{n,1}$, where $n \geq 2$.

Whether n is odd or even, there is
$$u_2 + u_N = q_{n-1} + q_{n-2} = q_{n,1}.$$

So, given $\alpha$ and $n$, there exists a unique $N \in (q_{n-1}, q_{n,1}]$, such that $N = u_2 + u_N$. In other words, in the interval $(q_{n-1}, q_{n,1}]$, there is and only one $N$ such that the cut to the Nth cut is of two gaps.

So, given $\alpha$ and $n$, when $q_{n-1} < N \le q_{n,1}$ (where $n \ge 2$), the number of times two gaps occur is at most 1.

Scenario 3: $N \le q_1$, then $u_j = j - 1$, $j = 1, 2, \cdots, N$.

In this case, the corresponding cake model would be the first loop of the cut, and the corresponding division of the interval $[0,1]$ would be the first layer truncation. In fact, in the range of values $1 < N \le q_1$, the gap cases of the cuts are all two.

So, given $\alpha$ and $n$, when $N \le q_1$, the number of times two gaps occur is at most $(q_1 - 1)$ times.

Let the frequency of occurrence of the two gaps be $p$, therefore

(i) If $q_{n,i-1} < N \le q_{n,i}$, where $2 \le i \le a_n$ and $n \ge 2$.

Then
$$p \le \frac{(q_1 - 1) + 1 + (a_2 - 1) + 1 + (a_3 - 1) + \cdots + 1 + (a_{n-1} - 1) + 1 + (a_n - 1)}{N}$$
$$= \frac{(a_1 + a_2 + \cdots + a_n) - 1}{N}$$
$$\in \left[ \frac{(a_1 + a_2 + \cdots + a_n) - 1}{q_{n,i}}, \frac{(a_1 + a_2 + \cdots + a_n) - 1}{q_{n,i-1}} \right)$$

(ii) If $q_{n-1} < N \le q_{n,1}$.

Then
$$p \le \frac{(q_1 - 1) + 1 + (a_2 - 1) + 1 + (a_3 - 1) + \cdots + 1 + (a_{n-1} - 1) + 1}{N}$$
$$= \frac{a_1 + a_2 + \cdots + a_{n-1}}{N}$$
$$\in \left[ \frac{a_1 + a_2 + \cdots + a_{n-1}}{q_{n,1}}, \frac{a_1 + a_2 + \cdots + a_{n-1}}{q_{n-1}} \right)$$

Since $q_{n,i-1} = (i-1)q_{n-1} + q_{n-2} > q_{n-1}$, where $i \ge 2$, then
$$\frac{(a_1 + a_2 + \cdots + a_n) - 1}{q_{n,i-1}} < \frac{a_1 + a_2 + \cdots + a_n}{q_{n-1}},$$

and
$$\frac{a_1 + a_2 + \cdots + a_{n-1}}{q_{n-1}} < \frac{a_1 + a_2 + \cdots + a_n}{q_{n-1}}.$$

Therefore, Theorem 1 holds by Lemma 6.

## 3  One example

Theorem 1 tells us that for almost all $\alpha \in (0,1)\setminus \mathbb{Q}$ in the Lebesgue velocimetry sense, the frequency of two gap lengths is 0 when $N \to \infty$, i.e., the frequency of three gap lengths is 1. In this regard, we find a specific class of $\alpha$ such that the frequency of two gaps is 0. This deepens our understanding.

**Theorem 7.**[6] **(Lagrang theorem)** Let $\alpha \in \mathbb{R}^+ \setminus \mathbb{Q}$, then the continued fractional expansion of $\alpha$ is periodic, if and only if, $\alpha$ is a root of a quadratic equation with integer coefficients.

**Theorem 8.** If $\alpha \in (0,1)\setminus \mathbb{Q}$ is a root of a quadratic equation with integer coefficients, then

$$\lim_{N \to \infty} \frac{\#\{1 \leq n \leq N \mid \#\Delta_n(\alpha) = 2\}}{N} = 0.$$

**Proof.** By theorem 7, let

$$\alpha = [0; b_1, b_2, \cdots\cdots, b_r, \dot{a}_1, a_2, \cdots\cdots, \dot{a}_k],$$

where

$$b_1, b_2, \cdots\cdots, b_r \in \mathbb{Z}_+ \cup \{0\}; \quad a_1, a_2, \cdots\cdots, a_k \in \mathbb{Z}_+$$

therefore

$$\begin{pmatrix} p_n & p_{n-1} \\ q_n & q_{n-1} \end{pmatrix}$$

$$= \begin{pmatrix} 0 & 1 \\ 1 & 0 \end{pmatrix} \begin{pmatrix} b_1 & 1 \\ 1 & 0 \end{pmatrix} \cdots \begin{pmatrix} b_r & 1 \\ 1 & 0 \end{pmatrix} \begin{pmatrix} a_1 & 1 \\ 1 & 0 \end{pmatrix} \cdots \begin{pmatrix} a_k & 1 \\ 1 & 0 \end{pmatrix} \begin{pmatrix} a_1 & 1 \\ 1 & 0 \end{pmatrix} \cdots \begin{pmatrix} a_k & 1 \\ 1 & 0 \end{pmatrix} \begin{pmatrix} a_1 & 1 \\ 1 & 0 \end{pmatrix} \cdots$$

let

$$n - r = jk + l,$$

Where

$$j \in \mathbb{Z}_+, 1 \leq l < k \text{ and } l \in \mathbb{Z}_+ \cup \{0\}.$$

then

$$\begin{pmatrix} p_n & p_{n-1} \\ q_n & q_{n-1} \end{pmatrix}$$

$$= \begin{pmatrix} 0 & 1 \\ 1 & 0 \end{pmatrix}\begin{pmatrix} b_1 & 1 \\ 1 & 0 \end{pmatrix}\cdots\begin{pmatrix} b_r & 1 \\ 1 & 0 \end{pmatrix}\left(\begin{pmatrix} a_1 & 1 \\ 1 & 0 \end{pmatrix}\cdots\begin{pmatrix} a_k & 1 \\ 1 & 0 \end{pmatrix}\right)^j \begin{pmatrix} a_1 & 1 \\ 1 & 0 \end{pmatrix}\begin{pmatrix} a_2 & 1 \\ 1 & 0 \end{pmatrix}\cdots\begin{pmatrix} a_l & 1 \\ 1 & 0 \end{pmatrix}$$

let

$$\begin{pmatrix} b_1 & 1 \\ 1 & 0 \end{pmatrix}\cdots\begin{pmatrix} b_r & 1 \\ 1 & 0 \end{pmatrix} = \begin{pmatrix} A & B \\ K & F \end{pmatrix}, \text{ where } A,B,K,F \in \mathbb{Z}_+ \text{ and } AF - BK = (-1)^r;$$

$$\begin{pmatrix} a_1 & 1 \\ 1 & 0 \end{pmatrix}\cdots\begin{pmatrix} a_k & 1 \\ 1 & 0 \end{pmatrix} = \begin{pmatrix} a & b \\ c & d \end{pmatrix} = T, \text{ where } a,b,c,d \in \mathbb{Z}_+ \text{ and } ad - bc = (-1)^k;$$

$$\begin{pmatrix} a_1 & 1 \\ 1 & 0 \end{pmatrix}\begin{pmatrix} a_2 & 1 \\ 1 & 0 \end{pmatrix}\cdots\begin{pmatrix} a_l & 1 \\ 1 & 0 \end{pmatrix} = \begin{pmatrix} e & f \\ g & h \end{pmatrix}, \text{ where } e,f,g,h \in \mathbb{Z}_+ \text{ and } eh - fg = (-1)^l.$$

therefore

$$\begin{pmatrix} p_n & p_{n-1} \\ q_n & q_{n-1} \end{pmatrix}$$
$$= \begin{pmatrix} 0 & 1 \\ 1 & 0 \end{pmatrix}\begin{pmatrix} A & B \\ K & F \end{pmatrix}\left(\begin{pmatrix} a & b \\ c & d \end{pmatrix}\right)^j \begin{pmatrix} e & f \\ g & h \end{pmatrix}.$$
$$= \begin{pmatrix} K & F \\ A & B \end{pmatrix} T^j \begin{pmatrix} e & f \\ g & h \end{pmatrix}$$

hence

$$|\lambda E - T| = \begin{vmatrix} \lambda - a & -b \\ -c & \lambda - d \end{vmatrix} = \lambda^2 - (d+a)\lambda + (ad - bc).$$

let

$$|\lambda E - T| = 0.$$

then

$$\lambda_1 = \frac{(d+a) + \sqrt{(d-a)^2 + 4bc}}{2},$$

$$\lambda_2 = \frac{(d+a) - \sqrt{(d-a)^2 + 4bc}}{2}.$$

evidently

$$\lambda_1 > 1,$$

$$\lambda_2 > 0 \Leftrightarrow ad - bc > 0 \Leftrightarrow k \text{ is even} \Leftrightarrow ad - bc = 1;$$

$$\lambda_2 < 0 \Leftrightarrow ad - bc < 0 \Leftrightarrow k \text{ is odd} \Leftrightarrow ad - bc = -1;$$

$$|\lambda_1| > |\lambda_2|.$$

Let $P$、$Q$ be second order matrices such that $T^j = \lambda_1^j P + \lambda_2^j Q \, (j \geq 0)$. So, let $j = 0, 1$, we have

$$\begin{cases} E = P + Q \\ T = \lambda_1 P + \lambda_2 Q \end{cases}.$$

Then

$$P = -\frac{1}{\sqrt{(d-a)^2 + 4bc}} \begin{pmatrix} \lambda_2 - a & -b \\ -c & \lambda_2 - d \end{pmatrix},$$

$$Q = \frac{1}{\sqrt{(d-a)^2 + 4bc}} \begin{pmatrix} \lambda_1 - a & -b \\ -c & \lambda_1 - d \end{pmatrix}.$$

Therefore

$$T^j = \frac{1}{\sqrt{(d-a)^2 + 4bc}} \begin{pmatrix} I & J \\ U & V \end{pmatrix} = \begin{pmatrix} u & v \\ w & z \end{pmatrix},$$

where

$$I = (ad - bc)\left(\lambda_2^{j-1} - \lambda_1^{j-1}\right) + a\left(\lambda_1^j - \lambda_2^j\right);$$

$$J = b\left(\lambda_1^j - \lambda_2^j\right);$$

$$U = c\left(\lambda_1^j - \lambda_2^j\right);$$

$$V = (ad - bc)\left(\lambda_2^{j-1} - \lambda_1^{j-1}\right) + d\left(\lambda_1^j - \lambda_2^j\right).$$

Hence

$$\begin{pmatrix} p_n & p_{n-1} \\ q_n & q_{n-1} \end{pmatrix} = \begin{pmatrix} K & F \\ A & B \end{pmatrix} \begin{pmatrix} u & v \\ w & z \end{pmatrix} \begin{pmatrix} e & f \\ g & h \end{pmatrix}.$$

Then

$$q_{n-1} = fAu + fBw + hAv + hBz$$

$$= \frac{1}{\sqrt{(d-a)^2 + 4bc}} \left( G\left(\lambda_1^{j-1} - \lambda_2^{j-1}\right) + M\left(\lambda_1^j - \lambda_2^j\right)\right),$$

where

$$G = (fA + hB)(bc - ad), \quad M = fAa + hBd + cfB + hAb.$$

Thus

$$\lim_{n \to \infty} \frac{a_1 + a_2 + \cdots + a_n}{q_{n-1}}$$

$$= \lim_{j \to \infty} \frac{a_1 + a_2 + \cdots + a_{jk+l+r}}{q_{jk+l+r-1}}$$

$$= \lim_{j \to \infty} \frac{(b_1 + b_2 + \cdots + b_r) + (a_1 + a_2 + \cdots + a_k) j + (a_1 + a_2 + \cdots + a_l)}{\frac{1}{\sqrt{(d-a)^2 + 4bc}} \left( G(\lambda_1^{j-1} - \lambda_2^{j-1}) + M(\lambda_1^j - \lambda_2^j) \right)}.$$

Let

$$(b_1 + b_2 + \cdots + b_r) + (a_1 + a_2 + \cdots + a_l) = R, \quad a_1 + a_2 + \cdots + a_k = S,$$

we have

$$\lim_{n \to \infty} \frac{a_1 + a_2 + \cdots + a_n}{q_{n-1}} = \lim_{j \to \infty} \frac{(R + Sj) \cdot \sqrt{(d-a)^2 + 4bc}}{G(\lambda_1^{j-1} - \lambda_2^{j-1}) + M(\lambda_1^j - \lambda_2^j)}.$$

By

$$\lambda_1^j - \lambda_2^j = (\lambda_1 - \lambda_2)(\lambda_1^{j-1} + \lambda_1^{j-2} \lambda_2 + \lambda_1^{j-3} \lambda_2^2 + \cdots + \lambda_1 \lambda_2^{j-2} + \lambda_2^{j-1}),$$

thus

$$\lim_{n \to \infty} \frac{a_1 + a_2 + \cdots + a_n}{q_{n-1}} \leq \lim_{j \to \infty} \frac{(R + Sj) \cdot \sqrt{(d-a)^2 + 4bc}}{G(\lambda_1 - \lambda_2) \lambda_1^{j-2} + M(\lambda_1 - \lambda_2) \lambda_1^{j-1}}$$

$$= \lim_{j \to \infty} \frac{R + Sj}{G \lambda_1^{j-2} + M \lambda_1^{j-1}}$$

$$= 0.$$

So

$$\lim_{n \to \infty} \frac{a_1 + a_2 + \cdots + a_n}{q_{n-1}} = 0.$$

i.e.

$$\lim_{N \to \infty} \frac{\#\{1 \leq n \leq N \mid \#\Delta_n(\alpha) = 2\}}{N} = 0.$$

□

# 4 Conclusion

This section builds on the work of Ravenstein[1] and works out the frequency problem of the three gap theorem by means of the related knowledge of continued fractions, leading to a result that is true almost everywhere, i.e., Theorem 1: For almost all $\alpha \in (0,1) \setminus \mathbb{Q}$ in the sense of the Lebesgue measure, there are

$$\lim_{N \to \infty} \frac{\#\{1 \leq n \leq N \mid \#\Delta_n(\alpha) = 2\}}{N} = 0.$$

For any $\alpha \in (0,1) \setminus \mathbb{Q}$, the frequency problem of the three-gap theorem remains worthy of continued study. We can further investigate the following problem: for any $\beta \in [0,1]$, let

$$Y = \left\{\alpha : \lim_{N \to \infty} \frac{\#\{1 \leq n \leq N \mid \#\Delta_n(\alpha) = 2\}}{N} = \beta\right\}.$$

What then is the Hausdorff dimension $\dim_H(Y)$ of the set $Y$? This problem would nicely relate the problem of the three gap theorem to fractal geometry and is therefore of some research interest.

[1]College of Mathematics and Statistics, Chongqing University, Chongqing, 401331, P. R. China
  Email address: zhx1327992788@163.com